\documentclass[10pt,reqno]{siamltex}

\usepackage{amssymb, amsmath, amsfonts, amscd}
\usepackage[english]{babel}
\usepackage{graphicx, epsfig}

\newcommand{\nm}{\noalign{\smallskip}}
\newcommand{\ds}{\displaystyle}

\newcommand{\p}{\partial}
\newcommand{\pd}[2]{\frac {\p #1}{\p #2}}
\newcommand{\eqnref}[1]{(\ref {#1})}

\newcommand{\Kbb}{\mathbb{K}}
\newcommand{\Rbb}{\mathbb{R}}
\newcommand{\Sbb}{\mathbb{S}}

\newcommand{\la}{\langle}
\newcommand{\ra}{\rangle}

\newcommand{\Ccal}{\mathcal{C}}
\newcommand{\Hcal}{\mathcal{H}}

\newcommand{\Kcal}{\mathcal{K}}

\newcommand{\Scal}{\mathcal{S}}

\newcommand{\Ga}{\alpha}
\newcommand{\Gb}{\beta}
\newcommand{\Gd}{\delta}
\newcommand{\Ge}{\epsilon}

\newcommand{\Gvf}{\varphi}

\newcommand{\Gl}{\lambda}

\newcommand{\Go}{\omega}

\newcommand{\GG}{\Gamma}

\newcommand{\GO}{\Omega}

\newcommand{\beq}{\begin{equation}}
\newcommand{\eeq}{\end{equation}}

\numberwithin{equation}{section}
\numberwithin{figure}{section}

\begin{document}
\title{Cloaking due to anomalous localized resonance in plasmonic structures of confocal ellipses\thanks{\footnotesize This work is supported by Korean Ministry of Education, Sciences and Technology through NRF grants Nos. 2010-0004091 and 2010-0017532.}}

\author{Daewon Chung\thanks{Department of Mathematics, Inha University, Incheon
402-751, Korea (chdaewon, hbkang, kskim, hdlee@inha.ac.kr).} \and Hyeonbae Kang\footnotemark[2] \and Kyoungsun Kim\footnotemark[2]  \and Hyundae Lee\footnotemark[2]}

\date{\today}
\maketitle

\begin{abstract}
If a core of dielectric material is coated by a plasmonic structure of negative dielectric material with non-zero loss parameter, then anomalous localized resonance may occur as the loss parameter tends to zero and the source outside the structure can be cloaked. It has been proved that the cloaking due to anomalous localized resonance (CALR) takes place for structures of concentric disks and the critical radius inside which the sources are cloaked has been computed. In this paper, it is proved that CALR takes place for structures of confocal ellipses and the critical elliptic radii are computed. The method of this paper uses the spectral analysis of the Neumann-Poincar\'e type operator associated with two interfaces (the boundaries of the core and the shell).
\end{abstract}



\section{Introduction}

We begin by defining the problem of the cloaking by anomalous localized resonance (CALR). Let $D$ and $\GO$ be bounded domains in $\Rbb^d$ $(d=2,3)$ such that $\overline{D} \subset \GO$, and for a parameter $\Gd>0$ let
\beq
\Ge_\Gd=\left\{\begin{array}{ll}
1& \textrm{in }\Rbb^d\setminus \overline{\GO}\,,\\
-1+i\Gd\quad&\textrm{in }\GO\setminus\overline{D}\,,\\
1& \textrm{in } D\,.
\end{array}\right.
\eeq
The number $\Gd$ is a loss parameter and $\Ge_\Gd$ represents the permittivity distribution. So, the configuration is regarded as a core $D$ with permittivity 1 coated by the shell $\GO\setminus\overline{D}$ with permittivity $-1+i\Gd$. For a given compactly supported function $f$ satisfying
\beq\label{fzerocond}
\int_{\Rbb^d} f \,dx=0\,,
\eeq
 we consider the following dielectric problem:
\begin{equation}\label{eqn:dproblem}
\nabla\cdot\Ge_\Gd\nabla V_\Gd= f\quad\textrm{in }\Rbb^2\,,
\end{equation}
with the decay condition $V_\Gd(x)\rightarrow 0$ as $|x|\rightarrow \infty$.  The condition \eqnref{fzerocond} represents conservation of charge and is required for the existence of the solution to \eqnref{eqn:dproblem} (with the decay condition). Let
\beq\label{power}
E_\Gd := \Im \int_{\Rbb^d} \Ge_\Gd |\nabla V_\Gd|^2\, dx = \int_{\GO\setminus D} \Gd |\nabla V_\Gd|^2\, dx
\eeq
($\Im$ for the imaginary part),  which approximately represents in quasistatic regime the time averaged
electromagnetic power produced by the source dissipated into heat (see \cite{acklm3}).

The problem of CALR is to characterize the source $f$ such that the following two conditions are satisfied:
\begin{equation}\label{eqn:blowup}
\lim_{\Gd \to 0} E_\Gd = \infty\quad\textrm{as }\Gd\rightarrow 0\,,
\end{equation}
and $V_\Gd/\sqrt{E_\Gd}$ goes to zero outside some
radius $a$, as $\Gd\to 0$:
\beq\label{bounded}
|V_\Gd (x)/\sqrt{E_\Gd}| \to 0  \quad\mbox{as } \Gd \to 0
\quad\mbox{when}\, |x| > a.
\eeq

A few words on the term CALR would be helpful. The condition \eqnref{eqn:blowup} is of spectral nature on which we discuss in the next section. It also reflects that the solution $V_\Gd$ oscillates very rapidly as $\Gd$ tends to zero. As shown in \cite{MN_PRSA_06, NMM_94}, the oscillation takes places near the interfaces $\p \GO$ and $\p D$. For this reason, this phenomenon is called the anomalous localized resonance. Its connection to the cloaking can be explained as follows. The condition \eqnref{eqn:blowup} implies that an infinite amount of energy dissipated
per unit time in the limit $\delta \to 0$ which is unphysical. So, if we scale the solution $V_\delta$ by the factor of $1/\sqrt{E_\delta}$ so that the energy dissipation becomes one independent of $\delta$, then according to \eqnref{bounded} the new solution approaches zero as $\Gd \to 0$. So, the effect of the source becomes negligible, in other words, the source is essentially invisible. We also say that the {\it weak} CALR occurs if
 \beq\label{weak_blowup1}
 \limsup_{\Gd \to 0} E_\Gd =\infty ,
 \eeq
and the limit in \eqnref{bounded} is replaced by limsup. It is worth mentioning that it is not known whether there is a source $f$ for which only the weak CALR, not CALR, takes place.

The phenomena of anomolous localized resonance was discovered by Nicorovici, McPhedran and
Milton \cite{NMM_94} and is related to invisibility cloaking
\cite{MN_PRSA_06}. It is also related to superlenses since, as shown in \cite{NMM_94}, the anomalous resonance
can create apparent point sources which look like mirror images.  For these connections and further developments on CALR, and for connections to cloaking due to complementary media, we refer to references cited in \cite{acklm3}.

It has been shown in a recent paper \cite{acklm1} that the condition \eqnref{eqn:blowup} (or \eqnref{weak_blowup1}) can be characterized in terms of the spectrum of the Neumann-Poincar\'e type (NP) operator associated with the double interfaces $\p D$ and $\p \GO$ of the equation \eqnref{eqn:dproblem}. This operator is a generalization of the classical NP operator associated with a single interface (see the next section). The spectral characterization is then applied to show, in particular, that if $D$ and $\GO$ are concentric disks in $\Rbb^2$
of radii $r_i$ and $r_e$, respectively, then there
is a critical radius
\beq
r_* = \sqrt{r_e^3/r_i}
\eeq
such that for any source $f$ supported outside $r_*$ CALR does not occur, and for sources $f$ supported inside weak CALR takes place. Furthermore, it is proved that if the Newtonian potential of the source satisfies an additional mild condition on its Fourier coefficients, then CALR takes place. The critical radius $r_*$ was also found in \cite{MN_PRSA_06}. These results were extended in
\cite{klsw} to the case when the core $D$ is not radial by a
different method based on a variational approach. There the source $f$ is assumed to be supported on circles.
The spectral approach was also applied to the case when the permittivity of the core is different from $1$ \cite{acklm2}. In the same paper it is proved that if $D$
and $\GO$ are concentric balls in $\Rbb^3$, CALR does not occur.
Such a discrepancy in two and three dimensions comes from the difference of the convergence rate
of the eigenvalues of the NP operators
associated with the structure. In 2D, they converge to $0$
exponentially fast, but in 3D they converge only at the rate of
$1/n$. The absence of CALR in such coated
sphere geometries is also linked with the absence of perfect plasmon waves:
see the appendix in \cite{klsw}. On the other hand, it is proved in \cite{acklm3} that it is possible to make CALR
occur in three dimensional spherical structure by using a shell with a specially
designed anisotropic dielectric constant. The constant is designed using the folded geometry. It is worth mentioning that
CALR is known to occur in a slab geometry in three dimensions \cite{MN_PRSA_06}. See also \cite{dong, Xiao}.

The circular structure seems the only known coated structure where CALR occurs, and it is of interest to find such a structure other than circular one. The purpose of this paper is to prove that if $D$ and $\GO$ are confocal ellipses, then CALR takes place. It is proved that like the circular case there is a critical elliptic radius (the elliptic radius is the radius in elliptic coordinates) such that for sources supported inside the radius CALR occurs, and it does not for sources supported outside the radius. We emphasize that the critical elliptic radius differs depending of elliptic radii of $D$ and $\GO$: if $\rho_i$ and $\rho_e$ are the elliptic radii of $D$ and $\GO$, respectively, then the critical radius is $(3\rho_e -\rho_i)/2$ if the shell is thin so that $\rho_e \le 3\rho_i$, and it is $2(\rho_e -\rho_i)$ if $\rho_e > 3\rho_i$. (See Theorem \ref{mainthm}.) This discrepancy comes from the difference of the asymptotic behavior of the eigenvalues of the NP operator. It is interesting to observe that if we put $\rho_e=\ln r_e$ and $\rho_i=\ln r_i$, then $(3\rho_e -\rho_i)/2 = \ln \sqrt{r_e^3/r_i}$, so the thin case is similar to the circular case. It is also interesting that the critical elliptic radius does not depend on the eccentricity of the ellipses, it only depends on the ratio of the elliptic radii.

To prove the main result of this paper we use the spectral approach developed in \cite{acklm1}. Based on computations of the NP operator on ellipses in \cite{AKL}, we are able to compute eigenvalues and eigenfunctions of the NP operator associated with confocal ellipses. We then prove Theorem \ref{mainthm} by adapting arguments similar to the circular case.

This paper is organized as follows. In section \ref{sec2}, we make a brief survey on spectral theory of the NP operator and review the spectral characterization of CALR in \cite{acklm1}. In section \ref{sec3} eigenvalues and eigenfunctions of the NP operator are computed. The main result and its proof are given in section \ref{sec4}.

\section{Neumann-Poincar\'e operator and a spectral characterization of CALR}\label{sec2}

In this section we review the result in \cite{acklm1} which characterizes the power $E_\Gd$ in terms of the spectral data of the Neumann-Poincar\'{e}-type operator.

Let $B$ be a bounded simply connected domain in $\Rbb^2$ with the Lipschitz boundary $\GG=\p B$. The single
layer potential $\mathcal{S}_\GG[\Gvf]$ of a function $\Gvf \in L^2(\GG)$ is defined by
\begin{align*}
\mathcal{S}_\GG[\Gvf] (x) &:= \frac{1}{2\pi} \int_{\GG} \ln |x-y| \Gvf (y) \, d
\sigma(y), \quad x \in \mathbb{R}^2.
\end{align*}
We also define a boundary integral operator $\Kcal_\GG^*$ on $L^2(\GG)$  by
$$
\Kcal_\GG^*[\Gvf](x) = \frac{1}{2\pi} \int_{\GG}
\frac{\la  x -y, \nu(x) \ra}{|x-y|^2} \Gvf (y)\,d\sigma(y), \quad x \in \GG.
 $$
Here and throughout this paper,
$\nu(x)$ is the outward unit normal  to $\GG$ at $x$, and $\la \;, \;\ra$ denotes the standard inner
product in $\mathbb{R}^2$. The operator $\Kcal_\GG^*$ is called the Neumann-Poincar\'{e} (NP) operator on $\GG$.

The NP operator is a fundamental object in solving classical Dirichlet (and Neumann) problems (see, for example, \cite{Folland76, KPS}). It is known that the spectrum of $\Kcal_\GG^*$ lies in $(-1/2, 1/2]$. If $\GG$ is $\Ccal^{1, \Ga}$ for some $\Ga>0$, then $\Kcal_\GG^*$ is compact, and hence the spectrum is discrete and accumulates at $0$. Recently there is revived interest in the spectrum of $\Kcal_\GG^*$ in connection with plasmonics. The dielectric constant $\Ge$ of the domain $B$ is called a plasmonic eigenvalue if $(\Ge+1)/2(\Ge-1)$ is an eigenvalue of $\Kcal_\GG^*$. We emphasize that $\Ge$ must be negative in order for $(\Ge+1)/2(\Ge-1)$ to lie in $(-1/2, 1/2]$, which reflects that $B$ is negative dielectric material. (See \cite{kang} and references therein.) It is worth mentioning that conventionally $2/\Gl$ is called a Fredholm eigenvalue if $\Gl$ is an eigenvalue of the NP operator.

The problem \eqnref{eqn:dproblem} involves two interfaces: $\GG_i:= \p D$ and $\GG_e:= \p\GO$.
Let $\Hcal=L^2(\GG_i)\times L^2(\GG_e)$ and $\Hcal_0=L^2_0(\GG_i)\times L^2_0(\GG_e)$ where $L^2_0(\GG_i)$ is the space of $L^2$-functions with the mean zero. The solution $V_\Gd$ to \eqnref{eqn:dproblem} can be represented as
\begin{equation}\label{solrep}
V_\Gd(x) = F(x) + \mathcal{S}_{\GG_i}[\Gvf_i](x) + \mathcal{S}_{\GG_e}[\Gvf_e](x),\quad x\in\mathbb{R}^2
\end{equation}
for a pair of potentials $(\Gvf_i, \Gvf_e) \in \Hcal_0$, where $F$ is the Newtonian potential of the source term $f$ in \eqnref{eqn:dproblem}, {\it i.e.},
 \begin{equation} \label{newton}
 F(x)= \frac{1}{2\pi} \int_{\Rbb^2} \ln |x-y| f(y) dy, \quad x \in \mathbb{R}^2.
 \end{equation}
It then follows from the continuity of the flux along interfaces that $(\Gvf_i, \Gvf_e)$ satisfies
\beq
(z_\Gd I + \Kbb^*) \begin{bmatrix} \Gvf_i \\ \Gvf_e \end{bmatrix} = g,
\eeq
where
\beq
z_\Gd  = \frac{i\Gd}{2(2-i\Gd)},  \quad
g=\begin{bmatrix}
\pd{F}{\nu_i} \\
\nm
-\pd{F}{\nu_e}
\end{bmatrix},
\eeq
and the operator $\Kbb^{\ast}$ is defined by
\beq\label{NPtype}
\quad \Kbb^{\ast}:=\begin{bmatrix}
\ds -\Kcal^{\ast}_{\GG_i} & \ds -\frac{\p}{\p \nu_i}\Scal_{\GG_e}\\
\ds \frac{\p}{\p \nu_e}\Scal_{\GG_i} & \ds \Kcal^{\ast}_{\GG_e}
\end{bmatrix} \, .
\eeq
Here and for the rest of this paper $\nu_i$ and $\nu_e$ stand for the outward normal to $\GG_i$ and $\GG_e$, respectively. This operator $\Kbb^{\ast}:\Hcal\to \Hcal$ is the Neumann-Poincar\'{e}-type operator associated with the interface problem \eqnref{eqn:dproblem} (with two interfaces $\GG_i$ and $\GG_e$).

Define the operator $\Sbb: \Hcal \to \Hcal$ by
\begin{equation} \label{eq S}
\Sbb = \begin{bmatrix}
\mathcal{S}_{\GG_i} & \mathcal{S}_{\GG_e} \\
\mathcal{S}_{\GG_i} & \mathcal{S}_{\GG_e}
\end{bmatrix}.
\end{equation}
It is helpful to make a comment on the operators on the off-diagonal in \eqnref{NPtype} and \eqnref{eq S}. For example, the operator $\Scal_{\GG_e}$ on the upper right corner in \eqnref{eq S} is an operator from $L^2(\GG_e)$ into $L^2(\GG_i)$. It is proven in \cite{acklm1} (see also \cite{kang}) that $-\Sbb$ is positive-definite on $\Hcal_0$ and $\Kbb^*$ is self-adjoint (and compact assuming that $\GG_i$ and $\GG_e$ are $\Ccal^{1,\Ga}$ for some $\Ga>0$) on $\Hcal_0$ equipped with the inner product
\beq
\la \Gvf, \psi \ra_\Sbb:= - \la \Gvf, \Sbb[\psi] \ra, \quad \Gvf, \ \psi \in \Hcal_0.
\eeq

Suppose $\ker\Kbb^{\ast} = \{0\}$ and let $\lambda_1,\lambda_2,\ldots$ ($|\lambda_1|
\ge |\lambda_2| \ge \ldots$) be the nonzero eigenvalues of $\Kbb^\ast$ and
$\Psi_n$ be the corresponding eigenfunctions normalized so that
$$
\la  \Psi_m,
\Psi_n \ra_\Sbb =\Gd_{mn}
$$
where $\Gd_{mn}$ is the Kronecker's delta. Then one can see easily that
\begin{equation}\label{phi_soluions}
\begin{bmatrix} \Gvf_i \\ \Gvf_e \end{bmatrix} =\ds \sum_n \frac{\la g, \Psi_n \ra_\Sbb}{\lambda_n + z_\Gd} \Psi_n.
\end{equation}
The following spectral characterization is proved in \cite{acklm1}:
\begin{equation}\label{gcloak}
E_\Gd\approx\Gd \sum_n \frac{|\la  g, \Psi_n \ra_\Sbb|^2}{\lambda_n^2 + \Gd^2} . \end{equation}
Here $a(\Gd) \approx b(\Gd)$ means that there exist positive constants $C_1$ and $C_2$ independent of $\Gd$ such that
$$
C_1 a(\Gd) \le b(\Gd) \le C_2 a(\Gd).
$$

\section{Spectrum of the NP operator on confocal ellipses}\label{sec3}

The elliptic coordinates $(\rho,\Go)$ for $x=(x_1, x_2)$ in the cartesian coordinates are defined by
\beq
x_1=R\cos\Go\cosh \rho, \ \ x_2=R\sin\Go\sinh \rho, \quad \rho> 0, \ 0\leq\Go\leq 2\pi\,.
\eeq
For any positive constant $\rho_0$ the set $E=\{(\rho,\Go): \rho =\rho_0\}$ is an ellipse whose foci are $(\pm R,0)$. So the confocal ellipses have the same elliptic coordinates. One can see easily that the length element and the normal derivative on $E$ are given by
 \begin{align}\label{length}
 d\sigma = \Xi \, d\Go \quad\mbox{and} \quad
 \frac{\p}{\p \nu} = \Xi^{-1}\frac{\p}{\p r},
 \end{align}
where
$$
\Xi= \Xi(\rho_0,\Go):= R\sqrt{\sinh^2 \rho_0 +\sin^2 \Go} \,.
$$

Let $E=\{(\rho,\Go): \rho =\rho_0\}$. It is proven in \cite{AKL} that if $h$ is a harmonic polynomial given by $h(x)=\cos n\Go(e^{n\rho}+e^{-n\rho})$ in elliptic coordinates for a non-negative integer $n$, then
\begin{equation}\label{eqn:SLP}
\Scal_{E}[\nabla h\cdot \nu](x)=\left\{\begin{array}{ll}
                                      (\Ga_n- \frac{1}{2})(e^{n\rho}+e^{-n\rho})\cos n\Go\,,\quad & \rho\leq \rho_0\,, \\
                                      \nm
                                      \Gb_n e^{-n\rho}\cos n\Go \,,& \rho>\rho_0\,,
                                    \end{array}
                                    \right.
\end{equation}
where
\beq\label{Gan}
\Ga_n=\frac{1}{2 e^{2n\rho_0}} \quad\textrm{and}\quad \Gb_n=\frac{-e^{2n\rho_0}+e^{-2n\rho_0}}{2}\,.
\eeq
Similarly, if $h(x)=\sin n\Go(e^{n\rho}-e^{-n\rho})\,,$ then
\begin{equation}\label{eqn:SLP2}
\Scal_{E}[\nabla h\cdot \nu](x)=\left\{\begin{array}{ll}
(-\Ga_n -\frac{1}{2})(e^{n\rho}-e^{-n\rho})\sin n\Go\,,\quad & \rho\leq \rho_0\,, \\
\nm
\Gb_n e^{-n\rho}\sin n\Go \,,& \rho>\rho_0\,.
\end{array}
\right.
\end{equation}
It then follows from the jump formula of the single layer potential, \eqnref{eqn:SLP} and \eqnref{eqn:SLP2} that
\begin{equation}\label{kstar_proj}
\Kcal^{\ast}_{E}[\Xi^{-1}\cos n\Go]=\Ga_n \Xi^{-1}\cos n\Go \quad\textrm{and}
\quad \Kcal^{\ast}_{E}[\Xi^{-1}\sin n\Go]=-\Ga_n \Xi^{-1}\sin n\Go.
\end{equation}
The jump formula of the single layer potential reads
$$
\frac{\p}{\p \nu}\Scal_{E}[\Gvf]\Big|_{-}(x)=\left(-\frac{1}{2}I+\Kcal^{\ast}_{E}\right)[\Gvf](x),\quad x\in E\,,
$$
where the subscript $-$ indicates the limit from the interior of $E$.

The formula \eqnref{kstar_proj} immediately yields the following lemma.
\begin{lemma}
$\Ga_n$ and $-\Ga_n$ $(n=0,1,2, \ldots)$ in \eqnref{Gan} are eigenvalues of the NP operator $\Kcal^{\ast}_{E}$ on the ellipse $E=\{(\rho,\Go):\rho=\rho_0\}$ and corresponding eigenfunctions are $\Xi(\rho_0,\Go)^{-1}\cos n\Go$ and $\Xi(\rho_0,\Go)^{-1}\sin n\Go$, respectively.
\end{lemma}

We now assume $\p D$ and $\p \GO$ are confocal ellipses whose common foci are denoted by $(\pm R,0)$.
Suppose that $\p D=\GG_i$ and $\p \GO=\GG_e$ are given by
$$
\GG_i=\{(\rho,\Go)\,:\,\rho=\rho_i\}\quad\textrm{and}\quad \GG_e=\{(\rho,\Go)\,:\,\rho=\rho_e\}\,.
$$
It is convenient to use the following notation: for $k=i, e$ and $n=0,1,2, \ldots$
\begin{equation}
\phi_n^{ck}(\Go):=\Xi(\rho_k,\Go)^{-1} \cos n\Go,\quad \phi_n^{sk}(\Go):=\Xi(\rho_k,\Go)^{-1} \sin n\Go.
\end{equation}
In view of \eqnref{length}, formula \eqnref{eqn:SLP} and \eqnref{eqn:SLP2} can be rewritten, respectively, as
\beq\label{SLP3}
\Scal_{\GG_k}[\phi_n^{ck}](x)=
\left\{\begin{array}{ll}
\ds - \frac{e^{n\rho}+e^{-n\rho}}{2n e^{n\rho_k}} \cos n\Go \,,\quad & \rho\leq \rho_k\,, \\
\nm
\ds -\frac{e^{n\rho_k}+ e^{-n\rho_k}}{2n e^{n\rho}} \cos n\Go \,,\quad & \rho > \rho_k\,,
\end{array}
\right.
\eeq
and
\beq\label{SLP4}
\Scal_{\GG_k}[\phi_n^{sk}](x)=
\left\{\begin{array}{ll}
\ds - \frac{e^{n\rho}-e^{-n\rho}}{2n e^{n\rho_k}} \sin n\Go \,,\quad & \rho\leq \rho_k\,, \\
\nm
\ds - \frac{e^{n\rho_k} - e^{-n\rho_k}}{2n e^{n\rho}} \sin n\Go \,,\quad & \rho > \rho_k\,,
\end{array}
\right.
\eeq
for $k=i,e$. So, we have
\begin{align*}
\frac{\p}{\p \nu_i}\Scal_{\GG_e}[\phi_n^{ce}]& =-\frac{e^{n\rho_i}-e^{-n\rho_{i}}}{2e^{n\rho_e}}\phi_n^{ci}, \\
\frac{\p}{\p \nu_e}\Scal_{\GG_i}[\phi_n^{ci}]& =\frac{e^{n\rho_i}+e^{-n\rho_i}}{2e^{n\rho_e}}\phi_n^{ce},\\
\frac{\p}{\p \nu_i}\Scal_{\GG_e}[\phi_n^{se}]& =-\frac{e^{n\rho_i}+e^{-n\rho_{i}}}{2e^{n\rho_e}}\phi_n^{si}, \\
\frac{\p}{\p \nu_e}\Scal_{\GG_i}[\phi_n^{si}]& =\frac{e^{n\rho_i}-e^{-n\rho_i}}{2e^{n\rho_e}}\phi_n^{se}.
\end{align*}
From these formula together with \eqnref{kstar_proj} we obtain
\begin{equation}\label{d:BNPO}
\Kbb^{\ast}\begin{bmatrix}
a\phi_n^{ci} \\
b\phi_n^{ce}
\end{bmatrix}=\begin{bmatrix}
\phi_n^{ci} & 0\\
0 & \phi_n^{ce}
\end{bmatrix} A
\begin{bmatrix}
a\\
b
\end{bmatrix},
\quad
\Kbb^{\ast}\begin{bmatrix}
a\phi_n^{si} \\
b\phi_n^{se}
\end{bmatrix}=\begin{bmatrix}
\phi_n^{si} & 0\\
0 & \phi_n^{se}
\end{bmatrix} B
\begin{bmatrix}
a\\
b
\end{bmatrix}\, ,
\end{equation}
where
$$
A= \begin{bmatrix}
\ds -\frac{1}{2e^{2n\rho_i}} & \ds\frac{e^{n\rho_i}-e^{-n\rho_i}}{2e^{n\rho_e}}\\
\ds \frac{e^{n\rho_i}+e^{-n\rho_i}}{2e^{n\rho_e}} & \ds \frac{1}{2e^{2n\rho_e}}
\end{bmatrix}, \quad
B=
\begin{bmatrix}
\ds \frac{1}{2e^{2n\rho_i}} & \ds \frac{e^{n\rho_i}+e^{-n\rho_i}}{2e^{n\rho_e}}\\
\ds \frac{e^{n\rho_i}-e^{-n\rho_i}}{2e^{n\rho_e}} & \ds -\frac{1}{2e^{2n\rho_e}}
\end{bmatrix} \, .
$$
Here $a$ and $b$ are constants. We first observe from \eqnref{d:BNPO} that since $A$ and $B$ are non-singular, $\ker\Kbb^{\ast} = \{0\}$. One can easily see that eigenvalues of the matrices $A$ and $B$ are those of $\Kbb^*$, and corresponding eigenfunctions can be computed using eigenvectors of $A$ and $B$. By computing eigenvalues and eigenvectors of $A$ and $B$, the following lemma is obtained.
\begin{lemma}\label{lem32}
The eigenvalues of $\Kbb^{\ast}$ are $\pm\Gl_{1,n}$ and $\pm\Gl_{2,n}$ $(n=0,1,2, \ldots)$ where
\begin{align*}
\Gl_{1,n} &= \frac{1}{4} \left( e^{-2n\rho_e}-e^{-2n\rho_i}-\sqrt{(e^{-2n\rho_e}-e^{-2n\rho_i})^2+4e^{-2n(\rho_e-\rho_i)}} \right),\\
\Gl_{2,n} &= \frac{1}{4} \left( e^{-2n\rho_e}-e^{-2n\rho_i}+\sqrt{(e^{-2n\rho_e}-e^{-2n\rho_i})^2+4e^{-2n(\rho_e-\rho_i)}} \right),
\end{align*}
and eigenfunctions (not normalized) corresponding to $\Gl_{1,n},-\Gl_{1,n},\Gl_{2,n},-\Gl_{2,n}$ are, respectively,
\begin{align}
\Psi_n^{1+}=\begin{bmatrix} a_{1,n}\phi_n^{ci} \\ b_n\phi_n^{ce}\end{bmatrix},\quad \Psi_n^{1-}=\begin{bmatrix}b_n \phi_n^{si} \\ a_{2,n}\phi_n^{se}\end{bmatrix},\quad
\Psi_n^{2+}=\begin{bmatrix} a_{2,n}\phi_n^{ci} \\ b_n\phi_n^{ce}\end{bmatrix},\quad \Psi_n^{2-}=\begin{bmatrix} b_n\phi_n^{si} \\ a_{1,n}\phi_n^{se}\end{bmatrix},\label{e:e_vector2}
\end{align}
where
\begin{align*}
a_{1,n}&=
{e^{-2n\rho_e}+e^{-2n\rho_i}+\sqrt{(e^{-2n\rho_e}-e^{-2n\rho_i})^2+4e^{-2n(\rho_e-\rho_i)}}} ,\\
a_{2,n}&=
{e^{-2n\rho_e}+e^{-2n\rho_i}-\sqrt{(e^{-2n\rho_e}-e^{-2n\rho_i})^2+4e^{-2n(\rho_e-\rho_i)}}}, \\
b_n &= -{2e^{-n(\rho_e-\rho_i)}(1+e^{-2n\rho_i})}.
\end{align*}
\end{lemma}

The formula \eqnref{eqn:SLP} and \eqnref{eqn:SLP2} show that
\begin{equation}
-\Sbb\begin{bmatrix}
a\phi_n^{ci} \\
b\phi_n^{ce}
\end{bmatrix}=\begin{bmatrix}
\ds \frac{\cosh n\rho_i}{ne^{n\rho_i}} & \ds \frac{\cosh n\rho_i}{ne^{n\rho_e}} \\
\ds \frac{\cosh n\rho_i}{ne^{n\rho_e}}  & \ds \frac{\cosh n\rho_e}{ne^{n\rho_e}}
\end{bmatrix}
\begin{bmatrix}
a\cos n\Go \\
b\cos n\Go
\end{bmatrix},
\label{d:BNPOS}\end{equation}
and
\begin{equation}
- \Sbb\begin{bmatrix}
a\phi_n^{si}\\
b\phi_n^{se}
\end{bmatrix}=\begin{bmatrix}
\ds \frac{\sinh n\rho_i}{ne^{n\rho_i}} & \ds \frac{\sinh n\rho_i}{ne^{n\rho_e}} \\
\ds \frac{\sinh n\rho_i}{ne^{n\rho_e}}  & \ds \frac{\sinh n\rho_e}{ne^{n\rho_e}}
\end{bmatrix}
\begin{bmatrix}
a\sin n\Go \\
b\sin n\Go
\end{bmatrix}.
\label{d:BNPOS2}\end{equation}
So, one can see from straight-forward calculations that $\Psi_n^{1\pm}, \Psi_n^{2\pm}$, $n=1,2,\ldots$, are orthogonal to each other with respect to the inner product $\la ~ ,\, \ra_\Sbb$, and
\begin{align*}
\la \Psi_n^{1+},\Psi_n^{1+}\ra_\Sbb  &= \frac{\pi}{n} \left( a_{1,n}^2 e^{-n\rho_i} \cosh n\rho_i + 2a_{1,n}b_n e^{-n\rho_e}\cosh n\rho_i+ b_n^2 e^{-n\rho_e}\cosh n\rho_e \right),\\
\la \Psi_n^{1-}, \Psi_n^{1-} \ra_\Sbb  &= \frac{\pi}{n} \left( b_n^2 e^{-n\rho_i} \sinh n \rho_i + 2a_{2,n}b_n e^{-n\rho_e}\sinh n\rho_i+ a_{2,n}^2 e^{-n\rho_e}\sinh n\rho_e \right),\\
\la \Psi_n^{2+},\Psi_n^{2+}\ra_\Sbb &= \frac{\pi}{n} \left( a_{2,n}^2 e^{-n\rho_i}\cosh n\rho_i + 2a_{2,n}b_ne^{-n\rho_e}\cosh n\rho_i+ b_n^2 e^{-n\rho_e}\cosh n\rho_e \right),\\
\la \Psi_n^{2-},\Psi_n^{2-}\ra_\Sbb &= \frac{\pi}{n} \left( b_n^2 e^{-n\rho_i}\sinh n\rho_i+2a_{1,n}b_ne^{-n\rho_e}\sinh n\rho_i + a_{1,n}^2 e^{-n\rho_e}\sinh n\rho_e \right).
\end{align*}

We now take a close look at the asymptotic behavior of eigenvalues as $n$ tends to $\infty$.
If $\GG_e$ and $\GG_i$ are sufficiently close to each other so that $\rho_e-\rho_i \le 2\rho_i$, then the asymptotic behavior of $\Gl_n$ is dominated by the term involved with $\rho_e-\rho_i$. More precisely, we have
$$
\sqrt{(e^{-2n\rho_e}-e^{-2n\rho_i})^2+4e^{-2n(\rho_e-\rho_i)}} = 2e^{-n(\rho_e-\rho_i)} + o(e^{-n(\rho_e-\rho_i)}).
$$
It then follows that
\begin{equation}\label{eigen_asymp1}
\Gl_{1,n}\sim -{e^{-n(\rho_e-\rho_i)}}\quad\textrm{and}\quad\Gl_{2,n}\sim {e^{-n(\rho_e-\rho_i)}}\,,
\end{equation}
where $p_n \sim q_n$ means that there are positive constants $C_1$, $C_2$, and $N$ such that
$$
C_1 q_n\le p_n \le C_2 q_n
$$
for all $n\ge N$. We also have
\begin{equation}
a_{1,n} \sim e^{-n(\rho_e-\rho_i)}, \quad a_{2,n} \sim -e^{-n(\rho_e-\rho_i)}, \quad b_n \sim -e^{-n(\rho_e-\rho_i)}.\label{ab_asymp1}
\end{equation}
On the other hand, if $2\rho_i<\rho_e-\rho_i$, we have
$$
\sqrt{(e^{-2n\rho_e}-e^{-2n\rho_i})^2+4e^{-2n(\rho_e-\rho_i)}} = e^{-2n\rho_i} +2 e^{-2n(\rho_e-2\rho_i)} + o(e^{-2n(\rho_e-2\rho_i)}),
$$
so that
\begin{equation}
\Gl_{1,n}\sim-{e^{-2n\rho_i}},\quad\Gl_{2,n}\sim{e^{-2n(\rho_e-2\rho_i)}}\,,\label{eigen_asymp2}
\end{equation}
and
\begin{equation}
a_{1,n} \sim e^{-2n\rho_i}, \quad a_{2,n}\sim -e^{-2n(\rho_e-2\rho_i)}, \quad b_n \sim -e^{-n(\rho_e-\rho_i)}.\label{ab_asymp2}
\end{equation}

Using \eqnref{ab_asymp1} and \eqnref{ab_asymp2}, we obtain the following lemma.
\begin{lemma}\label{lem_norm_est}
\begin{itemize}
\item[(i)] If $\rho_e \le 3\rho_i$, then
\begin{align}
\la \Psi_n^{1\pm}, \Psi_n^{1\pm}\ra_\Sbb,~ \la \Psi_n^{2\pm}, \Psi_n^{2\pm}\ra_\Sbb \sim  n^{-1} e^{-2n (\rho_e - \rho_i)}.\label{norm_est}
\end{align}
\item[(ii)] If $\rho_e > 3\rho_i$, then
\begin{align}
\la \Psi_n^{1+}, \Psi_n^{1+}\ra_\Sbb,~\la \Psi_n^{2-}, \Psi_n^{2-}\ra_\Sbb \sim  n^{-1} e^{-4n \rho_i}, \label{norm-asymp2}
\end{align}
and
\begin{align}
\la \Psi_n^{1-}, \Psi_n^{1-}\ra_\Sbb,~\la \Psi_n^{2+}, \Psi_n^{2+}\ra_\Sbb \sim  n^{-1} e^{-2n (\rho_e - \rho_i)}. \label{norm-asymp3}
\end{align}
\end{itemize}
\end{lemma}

\section{CALR on confocal ellipses}\label{sec4}

We now assume that the source $f$ is located outside $\GO$. Let us write the Newtonian potential of $f$ as
\begin{equation}
F(x)=c-\sum_{n\geq 1}(F_n^+\cos n\Go \cosh n \rho +F_n^-\sin n\Go \sinh n\rho)\, ,\label{e:NP1}
\end{equation}
where $c$ and $F_n^\pm$ are constants. Note that the series on the right-hand side converges in $0<\rho<\rho_0$ if and only if the coefficients $F_n^\pm$ satisfies
\begin{equation}\label{f_convergence}
\limsup_{n\to\infty} |F_n^\pm|^{1/n} \le e^{-\rho_0}.
\end{equation}
Since the source $f$ is located outside $\GO$ and hence the series converges in $0<\rho<\rho_e$, we infer that
\begin{equation}\label{Fpm_est}
\limsup_{n\to\infty} |F_n^\pm|^{1/n} \le e^{-\rho_e}.
\end{equation}

We now introduce a gap condition on $\{ F_n^\pm \}$ for confocal ellipses in analogy with the case of concentric disks or spheres in \cite{acklm1, acklm3, acklm2}. The sequence $F_n^\pm$ is said to satisfy the gap condition GC[$\rho_*$] for some constant $\rho_*$ if
\begin{quote}
GC[$\rho_*$]: there exists  a sequence $ \{n_k\} $ with
$n_1<n_2<\cdots$ such that
\beq
\lim_{k \to \infty} e^{-(n_{k+1}-n_k)(\rho_e -\rho_i)} e^{2{n_k}\rho_*}(|F_{n_k}^{+}|^2+|F_{n_k}^{-}|^2) = \infty .
\eeq
\end{quote}

The following is the main theorem of this paper.
\begin{theorem}\label{mainthm}
Let $f$ be the source function supported in $\Rbb^2 \setminus \overline{\GO}$ and $F$ be the Newtonian potential of $F$. Let
\beq\label{rhostar}
\rho_* =
\begin{cases}
\ds \frac{3\rho_e -\rho_i}{2} \quad &\mbox{if } \rho_e \le 3\rho_i, \\
2(\rho_e -\rho_i) \quad &\mbox{if } \rho_e > 3\rho_i.
\end{cases}
\eeq
\begin{itemize}
\item[(i)]
If $F$ does not extend as a harmonic function in $\{ \rho < \rho_* \}$, then
\beq\label{limsuped}
\limsup_{\Gd \to 0} E_\Gd=\infty
\eeq
and there is $C$ independent of $\Gd$ such that
\beq
|V_\Gd(x)| \le C
\eeq
for all $x$ satisfying $\rho \ge \rho_0$ provided that $\rho_0> 2\rho_e-\rho_i$ if $\rho_e \le 3\rho_i$ and $\rho_0> 3\rho_e-4\rho_i$ if $\rho_e < 3\rho_i$. So weak CALR takes place.
\item[(ii)] If, in addition, the
coefficients $F_n^\pm$ of $F$ in \eqnref{e:NP1}  satisfy $\mbox{GC}[\rho_{*}]$, then CALR takes place, {\it i.e.},
\beq
\lim_{\Gd \to 0} E_\Gd=\infty.
\eeq
\item[(iii)]  If $f$ is supported outside $\rho_*$ (so that $F$ extends as a harmonic function in $\{ \rho \le \rho_* \}$), then there is a constant $C$ such that
\beq
E_\Gd \le C
\eeq
for all $\Gd$.
\end{itemize}
\end{theorem}

Before proving Theorem \ref{mainthm} we make a few remarks on the condition GC[$\rho_*$]. As pointed out in \cite{acklm1}, it is a condition on the gap ($n_{k+1}-n_k$) among nonzero coefficients $F_n^\pm$. If nonzero coefficients are too sparse and the gap between them increases, GC[$\rho_*$] may fail. But the condition is satisfied by many source functions. For example, if $f$ is a dipole source, {\it i.e.}, $f(x)=a \cdot \nabla \Gd_{x_0}(x)$ for some $x_0$, then its Newtonian potential is given by
$$
F(x)= a \cdot \nabla_x G(x-x_0)
$$
where $G(x-x_0)= \frac{1}{2\pi}\ln|x-x_0|$. Using \eqnref{eqn:SLP} and \eqnref{eqn:SLP2}, one can see that $G(x-x_0)$ admits the expansion
\begin{equation}
G(x-x_0)=c- \sum_{n\ge 1} \frac{e^{-n\rho_0}}{n\pi} \left( \cos n\omega_0 \cos n\omega \cosh n\rho + \sin n\omega_0 \sin n\omega \sinh n\rho\right)
\end{equation}
for $x=(\rho,\omega)$ and $x_0=(\rho_0, \omega_0)$ with $\rho< \rho_0$. So $F$ satisfies GC[$\rho_*$] if $\rho_0 < \rho_*$, namely, the source $x_0$ is located inside $\rho_*$.

\medskip

\noindent{\sl Proof of Theorem \ref{mainthm}}. We only provide a proof for the case $\rho_e \le 3 \rho_i$. The other case can be proved without a major change.

Thanks to \eqnref{length}, we have
$$
\frac{\p F}{\p \nu_i}=\sum_{n\geq 1}(nF_n^+ \sinh{n\rho_i}\phi_n^{ci}+nF_n^- \cosh{n\rho_i}\phi_n^{si})
$$
and
$$
\frac{\p F}{\p \nu_e}=\sum_{n\geq 1}(nF_n^+ \sinh{n\rho_e}\phi_n^{ce}+nF_n^- \cosh{n\rho_e}\phi_n^{se}).
$$
So, $g$ is given by
$$
g:=\begin{bmatrix} \frac{\p F}{\p \nu_i}\\-\frac{\p F}{\p \nu_e}\end{bmatrix}=\sum_{n\geq 1}\begin{bmatrix} nF_n^+\sinh{n\rho_i}\phi_n^{ci}\\-nF_n^+\sinh{n\rho_e}\phi_n^{ce}\end{bmatrix}+\sum_{n\geq 1}\begin{bmatrix} nF_n^-\cosh{n\rho_i}\phi_n^{si}\\-nF_n^-\cosh{n\rho_e}\phi_n^{se}\end{bmatrix}\,.
$$
Therefore, we have
\begin{align}
\la g, \Psi_n^{1+}\ra_\Sbb
&=  \left( a_{1,n} e^{-n\rho_i}\cosh n\rho_i+b_ne^{-n\rho_e}\cosh n\rho_i\right)F_n^+\sinh n\rho_i\nonumber\\
 &\quad-\left(a_{1,n}e^{-n\rho_e}\cosh n\rho_i+b_n e^{-n\rho_e}\cosh n\rho_e \right)F_n^+\sinh n\rho_e, \label{g_proj1} \\
\la g, \Psi_n^{1-}\ra_\Sbb
 &=  \left( b_n e^{-n\rho_i}\sinh n\rho_i+a_{2,n}e^{-n\rho_e}\sinh n\rho_i\right)F_n^-\cosh n\rho_i\nonumber\\
 &\quad-\left(b_ne^{-n\rho_e}\sinh n\rho_i+a_{2,n} e^{-n\rho_e}\sinh n\rho_e \right)F_n^-\cosh n\rho_e,\label{g_proj2}\\
\la g, \Psi_n^{2+}\ra_\Sbb
&=  \left( a_{2,n} e^{-n\rho_i}\cosh n\rho_i+b_ne^{-n\rho_e}\cosh n\rho_i\right)F_n^+\sinh n\rho_i\nonumber\\
 &\quad-\left(a_{2,n}e^{-n\rho_e}\cosh n\rho_i+b_n e^{-n\rho_e}\cosh n\rho_e \right)F_n^+\sinh n\rho_e,\label{g_proj3}\\
\la g, \Psi_n^{2-}\ra_\Sbb
 &=  \left( b_n e^{-n\rho_i}\sinh n\rho_i+a_{1,n}e^{-n\rho_e}\sinh n\rho_i\right)F_n^-\cosh n\rho_i\nonumber\\
 &\quad-\left(b_ne^{-n\rho_e}\sinh n\rho_i+a_{1,n} e^{-n\rho_e}\sinh n\rho_e \right)F_n^-\cosh n\rho_e.\label{g_proj4}
\end{align}

With the aid of \eqnref{ab_asymp1},  the identities \eqnref{g_proj1}-\eqnref{g_proj4} yield
\begin{align}
&\la g, \Psi_n^{1+}\ra_\Sbb \sim  F_n^+ e^{n\rho_i},\quad \la g, \Psi_n^{1-}\ra_\Sbb \sim  F_n^- e^{n\rho_i},\nonumber\\
&\la g, \Psi_n^{2+}\ra_\Sbb \sim  F_n^+ e^{n\rho_i},\quad \la g, \Psi_n^{2-}\ra_\Sbb \sim  -F_n^- e^{n\rho_i}\label{proj_est}.
\end{align}
By \eqnref{phi_soluions} and Lemma \ref{lem32}, we have
$$
\begin{bmatrix} \Gvf_i \\ \Gvf_e \end{bmatrix} = \sum_{k=1,2}\sum_n
\left[ \frac{\la g, \Psi_n^{k+} \ra_\Sbb}{(\Gl_{k,n} + z_\Gd)\la \Psi_n^{k+}, \Psi_n^{k+} \ra_\Sbb} \Psi_n^{k+}
+ \frac{\la g, \Psi_n^{k-} \ra_\Sbb}{(-\Gl_{k,n} + z_\Gd) \la \Psi_n^{k-}, \Psi_n^{k-}\ra_\Sbb} \Psi_n^{k-} \right] \, .
$$
We thus have from \eqnref{e:e_vector2}
\begin{align}
|\mathcal{S}_{\GG_i}[\Gvf_i](x)| &\le \sum_{k=1,2}\sum_n
\left|\frac{a_{k,n} \la g, \Psi_n^{k+} \ra_\Sbb}{(\Gl_{k,n} + z_\Gd)\la \Psi_n^{k+}, \Psi_n^{k+} \ra_\Sbb} \mathcal{S}_{\GG_i}[\phi_n^{ci}](x) \right| \nonumber\\
& \quad +\sum_{k=1,2}\sum_n\left|\frac{b_{n} \la g, \Psi_n^{k-} \ra_\Sbb}{(-\Gl_{k,n} + z_\Gd) \la \Psi_n^{k-}, \Psi_n^{k-}\ra_\Sbb} \mathcal{S}_{\GG_i}[\phi_n^{si}](x) \right| \, . \label{solabs_est}
\end{align}
Note that
\begin{equation}
|\pm\Gl_{k,n} +z_\Gd| = \sqrt{\left(\pm\Gl_{k,n} - \frac{\Gd^2}{2(4+\Gd^2)}\right)^2 + \frac{\Gd^2}{(4+\Gd^2)^2} }\approx |\Gl_{k,n}|+\Gd
\label{lambdaplus}
\end{equation}
for all $n$. It thus follows from \eqnref{eigen_asymp1}, \eqnref{ab_asymp1}, \eqnref{norm_est} and \eqnref{proj_est} that
$$
|\mathcal{S}_{\GG_i}[\Gvf_i](x)| \le C\sum_n n e^{n(2\rho_e-\rho_i)} \left(\left|F_n^+\mathcal{S}_{\GG_i}[\phi_n^{ci}](x) \right|+\left|F_n^-\mathcal{S}_{\GG_i}[\phi_n^{si}](x) \right|\right).
$$
In view of \eqnref{SLP3} and \eqnref{SLP4}, we have
$$
|\Scal_{\GG_i}[\phi_n^{ci}](x)| + |\Scal_{\GG_i}[\phi_n^{si}](x)| \le C n^{-1} e^{n(\rho_i - \rho)}
$$
if $\rho > \rho_i$. So, we have
\begin{align}
|\mathcal{S}_{\GG_i}[\Gvf_i](x)|
& \le C\sum_n  e^{2n\rho_e}  e^{-n\rho}(|F_n^+|+|F_n^-|) \label{single_i_out}
\end{align}
for $\rho>\rho_i$.

Similarly we have, for $\rho>\rho_e$,
\begin{align}
|\mathcal{S}_{\GG_e}[\Gvf_e](x)|  &\le \sum_{k=1,2}\sum_n \left|\frac{b_n \la g, \Psi_n^{k+} \ra_\Sbb}{(\Gl_{k,n} + z_\Gd) \la \Psi_n^{k+}, \Psi_n^{k+}\ra_\Sbb} \mathcal{S}_{\GG_e}[\phi_n^{ce}] \right| \nonumber \\
& \quad + \sum_{k=1,2}\sum_n\left|\frac{a_{\hat{k},n} \la g, \Psi_n^{k-} \ra_\Sbb}{(-\Gl_{k,n} + z_\Gd)\la \Psi_n^{k-}, \Psi_n^{k-}\ra_\Sbb} \mathcal{S}_{\GG_e}[\phi_n^{se}] \right| \nonumber\\
&\le C\sum_n (|F_n^+|+|F_n^-|) e^{n(3\rho_e-\rho_i)} e^{-n\rho}. \label{single_e_out}
\end{align}
Here $\hat{k}=1$ if $k=2$ and $\hat{k}=2$ if $k=1$.

By \eqnref{solrep}, \eqnref{Fpm_est}, \eqnref{single_i_out}, and \eqnref{single_e_out} we conclude that
\begin{align*}
|V_\Gd(x)| \le |F(x)| + |\mathcal{S}_{\GG_i}[\Gvf_i](x)| + |\mathcal{S}_{\GG_e}[\Gvf_e](x)| < C
\end{align*}
regardless of $\Gd$ for any $\rho \ge \rho_0$ with $\rho_0 > 2\rho_e - \rho_i$.

Next we investigate the behavior of $E_\Gd$ as $\Gd\to 0$.  The formula \eqnref{gcloak} takes the form
\begin{align}
E_\Gd \approx \Gd \sum_{k=1,2}\sum_{n=1}^\infty \frac{1}{\Gl_{k,n}^2 + \Gd^2}\left(\frac{|\la  g,
\Psi_n^{k+} \ra_\Sbb|^2}{\la \Psi_n^{k+}, \Psi_n^{k+}\ra_\Sbb}+\frac{|\la  g,
\Psi_n^{k-} \ra_\Sbb|^2}{\la \Psi_n^{k-}, \Psi_n^{k-}\ra_\Sbb} \right).
\end{align}
Substituting \eqnref{eigen_asymp1}, \eqnref{norm_est}, and \eqnref{proj_est} into this formula yields
\begin{align}\label{Edelta}
E_\Gd \approx \sum_{n=1}^\infty \frac{\Gd n e^{2n\rho_e}((F_n^+)^2+(F_n^-)^2)}{e^{-2n(\rho_e -\rho_i)} + \Gd^2}.
\end{align}

Now we employ an argument similar to the one in \cite{acklm1} to determine if CALR takes place for given $f$.
Let $\rho_*$ be the number defined by \eqnref{rhostar} and suppose that $F$ does not extend as a harmonic function in $\{(\rho,\Go):\rho < \rho_*\}$. Then by \eqnref{f_convergence}, we have
\begin{equation*}
\limsup_{n\to\infty} |F_n^+|^{1/n} > e^{-\rho_*}\quad\mbox{or}\quad\limsup_{n\to\infty} |F_n^-|^{1/n} > e^{-\rho_*}.
\end{equation*}
So, there is a subsequence, say $\{ n_k \}$,  such that
\begin{equation}\label{limsup-cond}
e^{2 n_k \rho_*} ((F_{n_k}^+)^2+(F_{n_k}^-)^2) \ge 1
\end{equation}
for all $k$. Let $\Gd_k=e^{-n_k(\rho_e -\rho_i)}$ for $k=1,2,\ldots$. Then, we have
\begin{align*}
E_{\Gd_k} \approx\sum_{n=1}^\infty \frac{\Gd_k n e^{2n\rho_e}((F_n^+)^2+(F_n^-)^2)}{e^{-2n(\rho_e -\rho_i)} + \Gd^2}
& \ge   \frac{\Gd_k n_k e^{2n_k(2\rho_e-\rho_i)}}{2}((F_k^+)^2+(F_k^-)^2) \to \infty
\end{align*}
as $k \to \infty$. Thus \eqnref{limsuped} holds.

Suppose that $F_n^\pm$ satisfies GC[$\rho_*$] and $\{ n_k \}$ be the subsequence appearing in the condition. For $\Gd \to 0$, let $k(\Gd)$ be the number such that
$$
n_{k(\Gd)} \le - \frac{\ln \Gd}{\rho_e-\rho_i} < n_{k(\Gd)+1}.
$$
Then, we have
$$
\Gd > e^{-n_{k(\Gd)+1} (\rho_e-\rho_i)},
$$
and hence
 \begin{align*}
E_\Gd &\approx\sum_{n=1}^\infty \frac{\Gd n e^{2n\rho_e}((F_n^+)^2+(F_n^-)^2)}{e^{-2n(\rho_e -\rho_i)} + \Gd^2} \ge \frac{\Gd n_{k(\Gd)} e^{2n_{k(\Gd)}\rho_e}((F_{n_{k(\Gd)}}^+)^2+(F_{n_{k(\Gd)}}^-)^2)}{e^{-2{{n_k(\Gd)}}(\rho_e -\rho_i)}}\\
&\ge  n_{k(\Gd)} e^{-(n_{k(\Gd)+1}-n_{k(\Gd)})(\rho_e -\rho_i)} e^{2n_{k(\Gd)}\rho_*}(|F_{n_{k(\Gd)}}^{+}|^2+|F_{n_{k(\Gd)}}^{-}|^2)\rightarrow \infty\end{align*}
 as $\Gd\to 0$. This proves (ii).

 If the source $f$ is located outside $\rho_*$, then its Newtonian potential $F$ is harmonic in a neighborhood of $\{(\rho,\Go):\rho \le \rho_*\}$, and hence
$$
\limsup_{n\to\infty} |F_n^\pm|^{1/n} \le e^{-\rho_*-\Ge}
$$
for some $\Ge>0$. Thus it follows from \eqnref{Edelta} that
\begin{align*}
E_\Gd \le C   \sum_{n=1}^\infty { n e^{2n\rho_*}((F_n^+)^2+(F_n^-)^2)}<\infty.\end{align*}

This completes the proof.

\end{document}